\documentclass{amsart}
\usepackage{amssymb}
\usepackage{amsthm}
\usepackage{graphicx}
\newtheorem*{reftheorem}{Theorem}

\newtheorem{prop}{Proposition}

\numberwithin{equation}{section}
\title{Minimal number of self-intersections of the boundary of an
immersed surface in the plane}
\author{Larry Guth}
\address{Department of Mathematics, University of Toronto, 40 St. George St. Toronto ON M5S 2E4}
\email{lguth@math.toronto.edu}

\begin{document}
\begin{abstract} We find the minimal number of self-intersections
of the boundary of a surface of genus $g$ generically immersed in
$\mathbb{R}^2$.
\end{abstract}

\maketitle

Let $\Sigma$ be an oriented surface of genus $g \ge 1$ with one
boundary component.  We consider the class of all immersions $I:
\Sigma \rightarrow \mathbb{R}^2$ so that $I(\partial \Sigma)$
intersects itself transversely.  Among this class of immersions,
we determine the minimal number of self-intersections of
$I(\partial \Sigma)$.

\begin{prop} If $I$ is an immersion $I: \Sigma \rightarrow
\mathbb{R}^2$ and $I(\partial \Sigma)$ intersects itself
transversely, then $I(\partial \Sigma)$ has at least $2g+2$
self-intersections.  For each $g$, there is such an immersion so
$I(\partial \Sigma)$ has exactly $2g + 2$ self-intersections.
\end{prop}

This proposition answers a very 
simple case of a question that Gromov studied in the recent paper
\cite{G}.  Gromov gave estimates for the number of
self-intersections of the critical set of a generic map from one
manifold to another.  We can rewrite Proposition 1 in that
language as follows.  Suppose that $\Sigma'$ is a closed surface
of genus $2g$ without boundary.  Let $S$ be an embedded curve in
$\Sigma'$ which divides $\Sigma'$ into two surfaces each with
genus $g$.  It is possible to find a map $F$ from $\Sigma'$ to
$\mathbb{R}^2$ folded along the curve $S$ and with no other
singularities.  The curve $S$ is the singular set of the map $F$,
and $F(S) \subset
\mathbb{R}^2$ is the critical set of $F$.  Gromov observed that
as the topological complexity of $\Sigma$ increases, then the
topological complexity of the critical set $F(S)$ must also
increase.  As a corollary of Proposition 1, we see that for a
generic $F$ folded along $S$, the critical set must have at least
$2g + 2$ self-intersections, and this estimate is sharp.

\vskip5pt

{\bf Acknowledgements.} Yasha Eliashberg suggested this question
to me, and we had some helpful discussions about it.

\proof First we prove the lower bound.  The main ingredient of
the proof is the Whitney index formula, which relates the index
of an immersed curve with its self-intersections.  Whitney's
formula appears in his famous paper on immersed curves \cite{W}. 
(The more famous result of that paper is that any two immersed
curves with equal index are regular homotopic.)

Let $C$ be an oriented immersed curve given by an immersion
$\phi: S^1 \rightarrow \mathbb{R}^2$.  At any point $\theta$ of
$S^1$, the derivative of $\phi$ is a non-vanishing vector in
$\mathbb{R}^2$.  Therefore, the derivative of $\phi$ defines a
map from $S^1$ to $\mathbb{R}^2 -
\{ 0 \}$.  The winding number of this map is called the index of
the immersed curve.

We give $\mathbb{R}^2$ its standard orientation, and we orient
$\Sigma$ so that $I$ is orientation preserving.  We let $C$ be
the image $I(\partial \Sigma)$, with the boundary orientation. 
The first step of our proof is to show that the index of $C$ is
$1 -2g$.  This step follows from the Euler-Poincare formula.

Let $V$ be the pullback $I^*(\partial/\partial x)$. 
The vector field $V$ is a non-vanishing vector field on the
surface $\Sigma$.  We trivialize the bundle $T \Sigma$ restricted
to $\partial \Sigma$ so that the tangent vector to the boundary
is constant.  With respect to that trivialization, we let $W(V)$
be the winding number of the vector field $V$ along the boundary
$\partial \Sigma$.  According to the Euler-Poincare formula,
since $V$ is nowhere vanishing, $- W(V) = \chi(\Sigma) = 1 - 2g$. 
The immersion $I$ induces a trivialization of $T \Sigma$.  In
particular, it gives a second trivialization of $T \Sigma$ over
$\partial \Sigma$.  The index of $C$ is the
winding number of the tangent vector to $\partial \Sigma$ in this
second trivialization.  The second trivialization has the same
orientation as the first, and in the second trivialization, the
vector $V$ is constant.  Therefore, the winding number of the
tangent vector to $\partial \Sigma$ is $- W(V) = 1 - 2g$.

If $C$ is an oriented immersed curve with transverse
self-intersections, then its index and its self-intersections are
related by the Whitney index formula.  Let $p$ be a point of $C$
where the coordinate function $y$ achieves its minimum.  With
respect to $p$, we can give each self-intersection a sign $
\pm 1$.  If $x$ is a self-intersection, then at $x$ there are two
distinct unit tangent vectors tangent to the curve $C$ with the
correct orientation.  We call them $v_1$ and $v_2$.  We let $v_1$
be the tangent vector that occurs first if one follows the
immersed curve from the point $p$ until one reaches $x$. 
Finally, we say that $x$ is positive if $v_1$ is a positive
rotation from $v_2$.  The sign convention is illustrated in
Figure 1.

\includegraphics{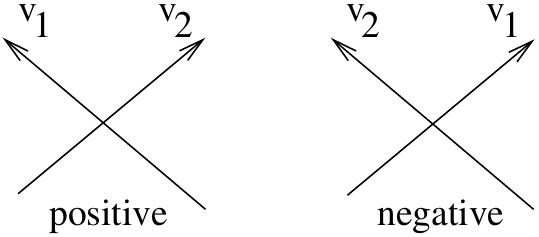}

We define $N^+$ to be the number of positive self-intersections
and $N^-$ to be the number of negative self-intersections. 
Finally, we define a number $\mu = \pm 1$ which depends on the
tangent vector of $C$ at $p$.  Because the function $y$ achieves
its minimum value at $p$, the tangent vector to $C$ at $p$ must
be $\pm \partial/\partial x$.  If the tangent vector is
$\partial/\partial x$, then $\mu = 1$, and if the tangent vector
is $- \partial / \partial x$, then $\mu = -1$.  In terms of these
conventions, the Whitney index formula reads as follows.

\begin{reftheorem}(Whitney) $ind(C) = \mu + N^+ - N^-.$
\end{reftheorem}

Figure 2 gives an example to illustrate the conventions and the
formula.  (Incidentally, this example bounds a surface of genus
1.)  For the curve in the figure, we have $ind(C) = - 1$, $\mu =
1$, $N^+ =1$, and $N^- = 3$.

\includegraphics{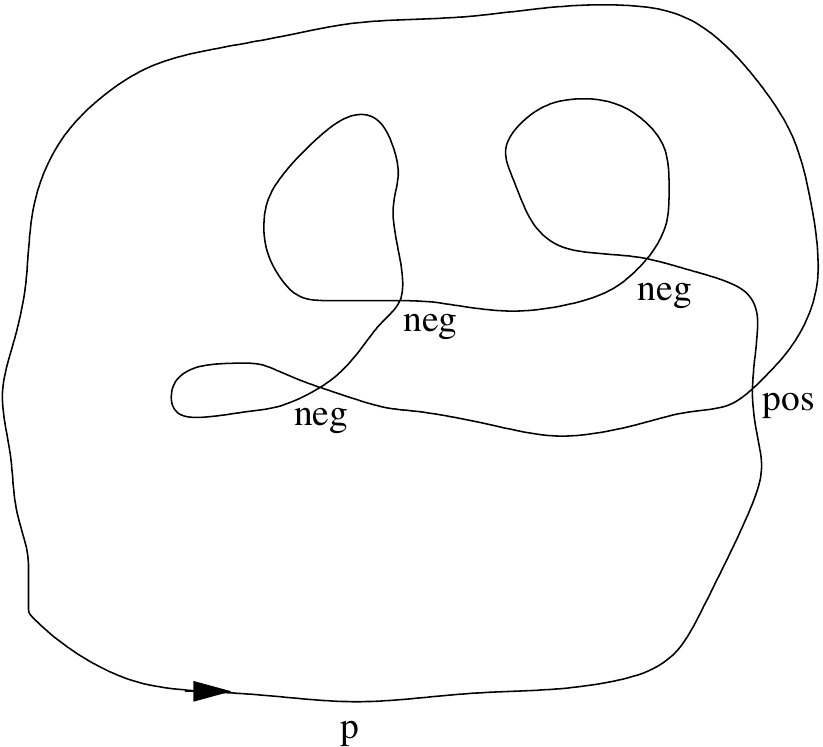}

For a general curve $C$ of index $1- 2g$, the Whitney index
formula shows that the total number of self-intersections, $N^- +
N^+$, is at least $2g - 2$.  Knowing only the index of $C$, this
estimate is the best possible, but using the immersed surface, we
can improve it in two places.

By translation, we can assume that the minimal value of $y$ on
$C$ is 0.  Since $I$ is an immersion, the minimal value of $y$ on
$I(\Sigma)$ occurs on $I(\partial \Sigma)$, and therefore the
image $I(\Sigma)$ lies above the line $y=0$.  Therefore, the
inward normal vector to $I(\Sigma)$ at $p$ must point in the
positive y-direction, and this implies that $\mu = +1$.  This is
the first improvement.

Because the surface $\Sigma$ has genus $g \ge 1$, the curve
$I(\partial \Sigma)$ must have at least one self-intersection. 
Let $x$ be the first point of self-intersection that one reaches
following the curve $I(\partial \Sigma)$ from $p$.  We claim that
the self-intersection at $x$ is positive.  This is the second
improvement.  Let $C_1$ denote the arc from $p$ to $x$, and let
$C_2$ denote a short piece of the other arc of $C$ through $x$. 
The positivity of the self-intersection is equivalent to knowing
that the inward normal vector of $I(\Sigma)$ along $C_2$ points
on the opposite side of $C_2$ from $C_1$.  But if the inward
normal vector lay on the same side as $C_1$, there would be a
second sheet of $I(\Sigma)$ under $C_1$, which would run down to
$p$ and then down past the line $y=0$, giving a contradiction. 
Therefore $N^+ \ge 1$.

According to the Whitney index formula, $N^- = \mu + N^+ -
ind(C) \ge 1 + 1 + (2g-1) = 2g + 1$.  Since we already showed
that $N^+ \ge 1$, the total number of self-intersections, $N^- +
N^+$, is at least $2g + 2$.  This finishes the proof of the lower
bound.

Next we construct an immersion with $2g + 2$ self-intersections,
for any $g$.  Our construction involves a few steps, illustrated
in Figure 3 in the case $g=2$.  We start with an immersion of the
disk with 2 self-intersections, illustrated in Figure 3.  Next,
we cut $g$ disjoint disks out of this immersed disk, removing
them from the multiplicity 1 region.  The result is an immersed
curve with $g+1$ components, bounding an immersed surface with
genus 0 and $g + 1$ boundary components.  The last step is to do
surgery on this immersed surface.  We glue in $g$ strips, as
shown in the figure.  Each strip connects one of the new circles
to the boundary of the original immersed disk.  The result is an
immersed circle with $2g + 2$ transverse self-intersections,
bounding an immersed surface of genus $g$.

\includegraphics{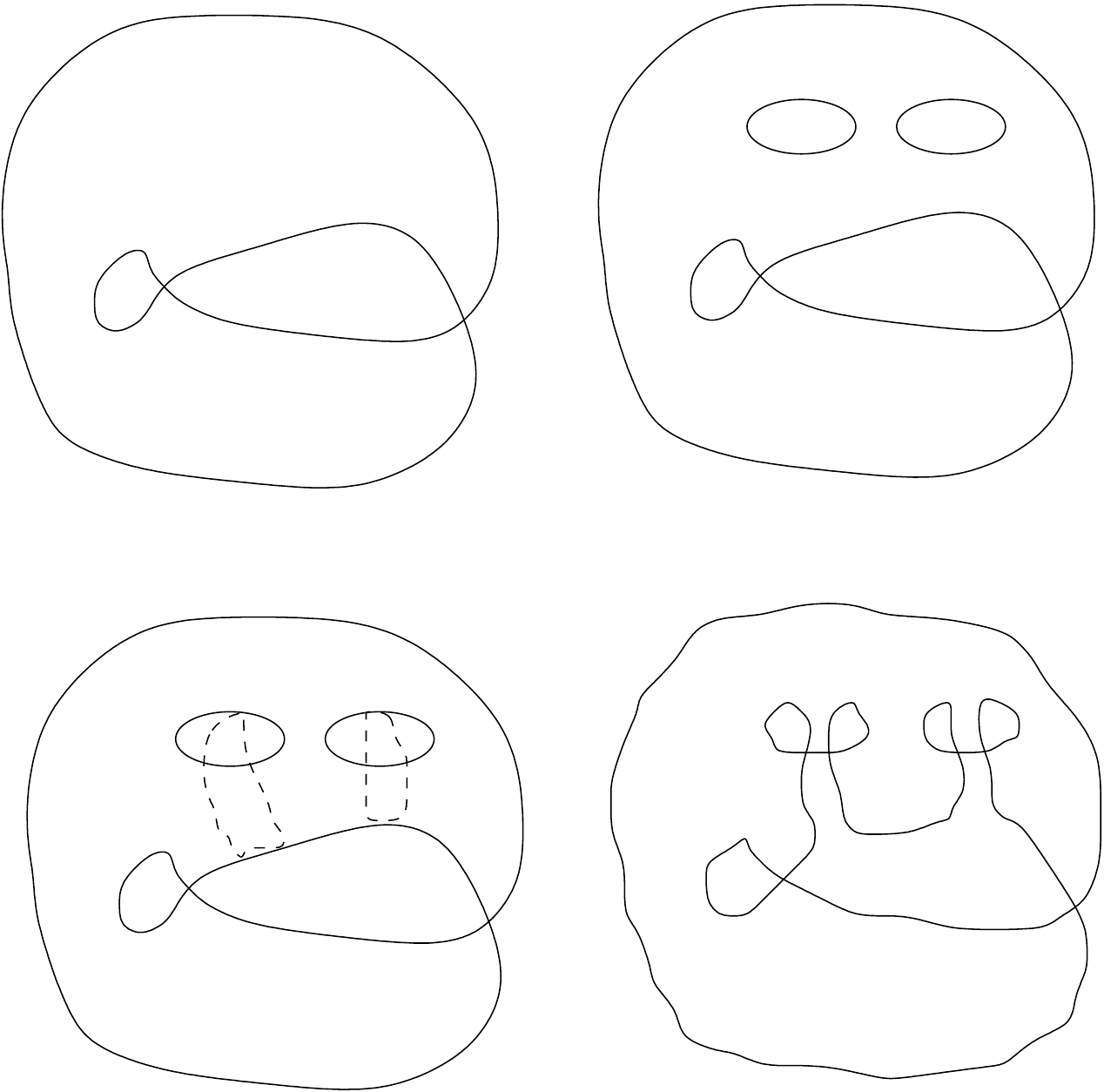}

\end{document}